\documentclass[12pt]{amsart}
\usepackage{amsmath,amscd,amssymb, amsthm,epsfig,psfrag}

\newtheorem{theorem}{Theorem}[section]
\newtheorem{lemma}[theorem]{Lemma}

\newtheorem{cor}[theorem]{Corollary}
\theoremstyle{definition}

\newtheorem{example}[theorem]{Example}
\newtheorem{claim}[theorem]{Claim}
\newtheorem{conjecture}[theorem]{Conjecture}

\def\H3{$\mathbb H^3$}
\def\n{$\nu$}

\def\L{$\mathcal{L}$\ }

\newcommand{\RR}{\mathbb{R}}

\newcommand{\QQ}{{\mathbb Q}}
\newcommand{\PP}{{\mathbb P}}
\newcommand{\interior}{{\rm int\ }}

\newcommand{\ZZ}{{\mathbb Z}}
\newcommand{\neigh}{{\mathcal N}}

\begin{document}

\title{Lower bounds on volumes of hyperbolic Haken 3-manifolds} 

\author[Ian Agol]{%
        Ian Agol} 
\address{%
        Department of Mathematics \\
        U. C. Davis \\
        Davis, CA 95616} 
\email{%
        iagol@math.ucdavis.edu}  
%




 

\maketitle
\section{Introduction}
In this paper, we find lower bounds for the volumes of certain hyperbolic 
Haken 3-manifolds.  The 
theory of Jorgensen and Thurston shows that the volumes of hyperbolic 3-manifolds 
are well-ordered, but no one has been able to find the smallest one. 
 The best known result for closed manifolds is that the smallest 
closed hyperbolic 3-manifold has volume $>0.16668$, proven by Gabai, 
Meyerhoff, and Thurston \cite{GMT}. Their proof involves extensive rigorous 
computations. The smallest  known closed 3-manifold is the Weeks manifold $W$, 
$Vol(W) \approx 0.94270$ \cite{W}, which is the 2-fold branched cover
over the knot $9_{49}$ \cite{Ve}.

One could also ask for the smallest hyperbolic manifolds with certain 
characteristics. The smallest non-compact, disoriented hyperbolic 3-manifold 
is the Gieseking manifold, which is double covered by the figure eight knot 
complement, proven by Colin Adams \cite{A}. Its volume = V$_3\approx 1.01494$,
 since it is obtained by pairwise gluing the faces of the regular ideal 
tetrahedron in $\mathbb H^3$, which has this volume. A recent result of Cao 
and Meyerhoff \cite{CM} shows that the smallest oriented noncompact hyperbolic
 3-manifolds are the figure-eight knot complement and its sister manifold, 
which have volume $\rm 2V_3\approx 2.03$. Kojima and Miyamoto \cite{KM,Mi}
have found the smallest hyperbolic 3-manifolds with totally geodesic
 boundary, which include Thurston's
tripos manifold \cite{Th}. Culler and Shalen have a series of papers deriving 
lower bounds for volumes of closed hyperbolic 3-manifolds M, where 
$\dim(H_1(M;{\mathbb Q}))=\beta_1(M) \geq1$. They found $Vol(M)>.34$, 
where one
 excludes ``fibroids'' if $\rm \beta_1(M)=1$ \cite{CS2,CS1}.
Culler, Shalen, and Hersonsky showed that if $\beta_1(M)\geq 3$, 
then $Vol(M)\geq .94689 > Vol(W)$, which shows that the smallest
volume $3$-manifold must have $\beta_1(M)\leq 2$ \cite{CHS}. Joel Hass has shown
that there is an upper bound to the genus of an acylindrical surface
in terms of the volume of a manifold \cite{H}. One consequence
of our results is that if $\beta_1(M)=2$ or $\beta_1(M)=1$ and
$M$ is not fibred over $S^1$, then $Vol(M)\geq \frac45 V_3$. 
If $M$ contains an acylindrical surface $S$, then $Vol(M)\geq -2V_3\chi(S)$.

Section  2 gives the necessary definitions and the
statement of the main theorem. Section 3 gives some examples.
Sections 4-6 are devoted to the proof of the main theorem.
Sections 7-9 give applications of the main theorem. 

Acknowledgements: We thank Daryl Cooper, Mike Freedman, Darren Long, Pat
Shanahan, and Bill Thurston for helpful conversations.

\section{Definitions and statement of the main theorem}
We will assume that all $3$-manifolds are orientable, and
usually will be homeomorphic to a compact manifold with
some subsurface of the boundary removed. A properly embedded 
{\it incompressible} surface $S$ in $M^3$ is a two-sided
surface for which the fundamental group injects (excluding
$S^2$). A manifold
is {\it irreducible} if every embedded $S^2$ bounds a ball. An irreducible
manifold with an incompressible surface is called {\it Haken}.  It is 
{\it hyperbolic} if it has a 
complete riemannian metric of constant sectional curvature -1, or 
equivalently, its universal cover is isometric to hyperbolic space \H3. 
If the manifold has boundary, then it is hyperbolic with totally
geodesic boundary if the metric is locally modelled on a close half-space 
in \H3 bounded by a geodesic plane.

The Thurston norm of a connected surface $S$ is $\max\{-\chi(S),0\}$,
and extends to disconnected surfaces by summing over the 
components. An oriented surface which represents a non-trivial
homology class in $H_2(M;\QQ)$ is Thurston norm-minimizing if
it is has minimal Thurston norm over all surfaces in the 
same homology class. The norm on the sublattice represented
by embedded surfaces extends by linearity to all of 
$H_2(M;\QQ)$.  Thurston showed that for a hyperbolic 
3-manifold, this norm on homology is indeed a norm \cite{Th3}.
Also, if $\beta_2(M)\geq 2$, he showed that there is 
some homology class with a surface which is not a fiber
over $S^1$. 

$M$ is {\it atoroidal} if any $\pi_1$-injective mapping of
a torus into $M$ is homotopic into $\partial M$. A {\it pared
acylindrical manifold} is a pair $(M,P)$ where
\begin{itemize}
\item
$M$ is a compact, irreducible, atoroidal $3$-manifold,
\item
$P\subset \partial M$ is a union of incompressible annuli
and tori, such that every map 
$(S^1\times I,S^1\times\partial I)\to(M,\partial M-P)$ that is $\pi_1$-injective deforms as a map of pairs into $P$.
\end{itemize}
$P$ is called the parabolic locus of the pared manifold $(M,P)$.
Let $\partial_0 M=\partial M- P$.
A theorem of Thurston shows that if $(M,P)$ is pared 
acylindrical, and $M$ is not a torus or a ball, 
then $M\setminus P$ has a hyperbolic metric with totally
geodesic boundary $\partial_0 M$ \cite{Mo}.

An $I-bundle$ pair is a pair $(I-bundle,\partial I-bundle)$
over a surface. Let $M$ be an irreducible $3$-manifold, and
let $Q\subset \partial M$ be an incompressible surface. There is a
subpair $(\Sigma, S)\subset (M,Q)$ called the {\it characteristic
sub-pair} of $(M,Q)$ \cite{Jac,Mo}.
It is uniquely determined up to isotopy of pairs by the following conditions:
\begin{itemize}
\item
Each component of $(\Sigma,S)$ either is an $(I-bundle,\partial I-bundle)$-pair
or has a Seifert fibration structure in which $S$ is a union of fibers in
the boundary.
\item
The components of $\partial_1\Sigma = \partial\Sigma \setminus \interior S$
are essential annuli and tori in $(M,Q)$. Each is either a component of
$\partial M\setminus \interior Q$ or is not parallel into $\partial M$.
\item
No component of $(\Sigma,S)$ is homotopic in $(M,Q)$ into a
distinct component.
\item
Any map of the torus or the annulus into $(M,Q)$, which is injective
on the fundamental group and essential as a map of pairs, is
homotopic as a map of pairs into $(\Sigma, S)$. 
\end{itemize}

We will call the union of components consisting of $I-bundles$ the
{\it characteristic I-bundle}. 
 If $M$ is atoroidal, $P\subset \partial M$ consists of tori and annuli,
 and $(\Sigma,S)\subset (M,\overline{\partial_0 M})$ is the characteristic
sub-pair, then $Guts(M,P)=(\overline{M\setminus \Sigma},
(\overline{M\setminus\Sigma}) \cap \Sigma)$ is a pared
acylindrical manifold pair. If $(M,P)$ is already pared acylindrical, then
the characteristic sub-pair of $(M,\partial_0 M)$ is a product
neighborhood of $P$, and $Guts(M,P)$ is homeomorphic to $(M,P)$. 
If $N=M\setminus P$, for some $(M,P)$, then we define $Guts(N)=Guts(M,P)$,
{\it i.e.} the pared locus is implicit for such non-compact manifolds. 

\begin{theorem} \label{main}
If $M$ is a hyperbolic manifold containing an incompressible
surface $S$, then $Vol(M)\geq -2V_3\chi(Guts(M\setminus\neigh(S)))$. 
\end{theorem}

In particular, if $M\setminus\neigh(S)$ is acylindrical, then 
$Vol(M)\geq -2V_3\chi(S)$. 

We conjecture a sharper bound. If $\partial N$ is acylindrical,
let $Vol(N)$ be the volume of the components of
$N$ which have a hyperbolic metric with
totally geodesic boundary. 
\begin{conjecture} \label{conj1}
If $M$ is a hyperbolic manifold containing an incompressible
surface $S$, then $Vol(M)\geq Vol(Guts(M\setminus\neigh(S)))$.
\end{conjecture}

Let $V_{oct}\approx 3.66$ denote the volume of a regular
ideal octahedron in \H3. Miyamoto showed that if $M$ is hyperbolic 
with totally
geodesic boundary, then $Vol(M)\geq -V_{oct}\chi(M)$, with
equality holding only for manifolds made up of regular ideal
octahedra, glued together in the pattern of an ideal
trianglulation \cite{Mi}. Then conjecture \ref{conj1} and
Miyamoto's result  would imply that we could
improve the constant in theorem \ref{main} from 
$2V_3$ to $V_{oct}$.

\section{Examples}

\begin{example}
The manifold obtained by $(\frac{14}{3},\frac{3}{2})$ filling on the
Whitehead link is a Haken manifold of volume 2.2077. That this is
Haken was shown by Nathan Dunfield (personal communication). 
He also found that $+10$ filling
on the knot $5_2$ is Haken and has volume $2.3627$. These are the
smallest volume closed Haken 3-manifolds that I know of. Conjecture \ref{conj1}
would imply that these manifolds have no incompressible surfaces $S$
with $\chi(Guts(M\setminus\neigh(S)))<0$.

\begin{figure}[htb]
	\begin{center}
	\psfrag{a}{$\frac{14}{3}$}
	\psfrag{b}{$\frac32$}
	\psfrag{c}{$10$}
	\epsfxsize=3in
	\epsfbox{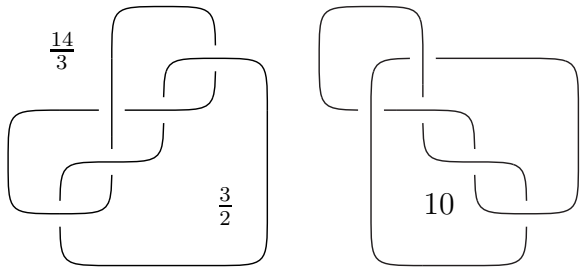}
	\caption{Dunfield's small volume hyperbolic Haken 3-manifolds}
	\end{center}
\end{figure}
\end{example}
\begin{example}
A manifold containing incompressible surfaces of unbounded genus can have 
bounded volume.
For example, if $\beta_1(M)\geq 2$, then any primitive $\ZZ$-homology class
is represented by a surface of minimal Thurston norm \cite{Th3}.  
Since the image of $H_2(M;\ZZ)$ in $H_2(M;\RR)$ is a lattice, 
there are primitive elements of unbounded norm,
and therefore corresponding surfaces of unbounded genus.
\end{example}
\begin{example}
In section \ref{planar} we show that if a hyperbolic link complement
contains an incompressible punctured sphere with meridian 
boundary components, then the volume is bounded below by $4V_3$. 
The smallest volume example that we know of (along with
its mutant) is given in figure
\ref{mutant}.  There is a unique ideal triangulation of $B^3$ with
one tetrahedron. The one-skeleton is a tangle in $B^3$. Viewing
a regular ideal octahedron as a truncated tetrahedron, we can
glue it together according to the  pattern of the
triangulation to obtain a geometric
structure on the tangle complement with boundary a geodesic
4-punctured sphere. The double gives the link
in figure \ref{mutant}. Conjecture \ref{conj1} would
imply that it and its mutant 
are the  minimal volume hyperbolic links containing a meridianal
incompressible planar surface.
 \begin{figure} [htb]
	\begin{center}
	\epsfxsize=2in
	\epsfbox{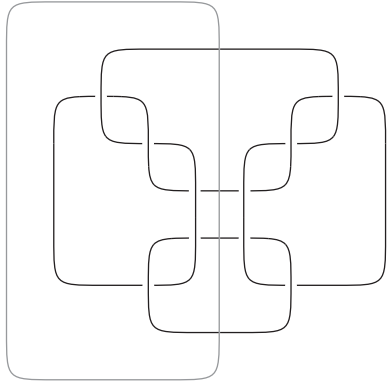}
	\caption{\label{mutant} A link of volume 7.32 with a Conway sphere}
	\end{center}
\end{figure}
\end{example}
\begin{example}\label{twobridge}
We can find lower bounds for volumes of certain 2-bridge    
links. A $2$-bridge link $L$ is obtained by taking the boundary of
a plumbing of $n$ twisted annuli and moebius strips 
(see figure \ref{2bridg}). If the number
of twists in each strip is $\geq 3$, then $Vol(S^3\setminus L)
\geq 2V_3(n-1)$. A more general bound will be given in
section \ref{sec2bridge}.
\end{example}
\begin{example}
If $M$ contains a surface $S$ which is Thurston norm-minimizing
and geometrically finite, then $Vol(M)\geq \frac45 V_3$. This
will be proven in section \ref{cyclic2}. This
occurs when $\beta_2(M) \geq 2$, or 
$\beta_1(M)=1$ and $M$ does not fiber over $S^1$. If conjecture \ref{conj1} is
true, then $Vol(M)\geq \frac25 V_{oct} > Vol(W)$, which would imply
that this type of 3-manifold could not be smallest volume. 
\end{example}

\section{Generalized Gromov norms}
Here, we discuss the needed properties of generalized Gromov
norms. Gromov introduced a norm as
a topological invariant of homology classes. In particular, ones
gets an invariant
of the fundamental class of a manifold \cite{Gr}. For hyperbolic manifolds,
he proved that the volume is proportional to the norm of the
fundamental class. We need a generalization of this to cell
cycles. Cell cycles are to simplicial cycles as cell 
complexes are to simplicial complexes.  Similar to 
the Gromov norm, one can take the norm on cycles represented
by particular polyhedral cells with simplicial faces. 
Then for hyperbolic 3-manifolds, the
volume will also be proportional to the norm of the fundamental
class, where the proportionality constant is just the volume of the
largest volume geodesic representitive of the polyhedron in \H3 .

For us, a cell will be a piecewise linear ball $P$ whose boundary
is a cell complex $\partial P$. Cells are given orientations in your favorite
fashion, {\it e.g.} by orienting the simplices of a triangulation
or by orienting the frame bundle. 
For an oriented cell $P$, we'll let $\overline{P}$
be the same cell with the opposite orientation, and 
if $\sigma: P\to M$ is a map, then $\overline{\sigma}:\overline{P}\to M$
is the same map on the polyhedron with opposite orientation.
Let $\mathcal P$ be a collection of cells, such that the faces
of any cell in $\mathcal P$ are also in $\mathcal P$.
Define the chain complex
 $$C_*(M;\mathcal{P},\RR)= \{\sum_{\sigma: P\to M} 
a_\sigma \sigma : a_\sigma\in \RR, P\in \mathcal{P} \}/\{
\sigma = -\overline{\sigma}\}.$$ Each $n-1$-cell in the cell complex 
$\partial P$ inherits an orientation from $P$. 
We define $\partial(\sigma)=\sum_{K\in \partial P} \sigma|_K$, and 
extend $\partial$ to $C_*(M;\mathcal{P},\RR)$ linearly. Then one can check that 
$(C_*(M;\mathcal{P},\RR), \partial)$ is a chain complex. The homology
of this complex might be different from singular homology. 

For example, 
we may take the cell complex $\mathcal O$ of the octahedron. More generally,
we can take any convex polyhedron $P$ in \H3 with triangular faces (the
dihedral angles between some faces may be $\pi$). Then a theorem 
of Thurston says that this polyhedron has a geometric realization in
\H3 as an ideal polyhedron of maximal volume, when one can place
horospheres centered at each vertex which are tangent along each edge
of the polyhedron (like a bunch of frog eggs). Let $V(P)$ be this maximal
volume, and let $\mathcal P$ be the corresponding cell complex of $P$. 
We can now take the Gromov norm with respect to $\mathcal P$. 
For cycles in $Z_*(M;\mathcal{P},\RR)$, we define the norm to be
$||\sum_{\sigma: P\to M} a_\sigma \sigma||=\sum_{\sigma: P\to M} |a_\sigma|$, where 
$P\in \mathcal{P}, a_\sigma\in \RR$.
Since $\mathcal{P}$ only has cells up to dimension 3, 
$H_3(M;\mathcal{P},\RR)=Z_3(M;\mathcal{P},\RR)$. But we have
a map $\psi:H_3(M;\mathcal{P},\RR)\to H_3(M;\RR)$ (singular homology),
given by choosing some subdivision of the polyhedron $P$ into 
oriented simplices (the map $\psi$ should average over all twelve
parametrizations of a simplex, so that we indeed get a singular 
cycle). For a cycle $z\in H_3(M;\mathcal{P},\RR)$, 
we will say $z$ is a fundamental cycle if $[\psi(z)]=[M]$. If
$z$ is transverse to a point in $M$, then the degree of $z$ at that
point will be 1. Now, we define 
$||[M]||_P= \inf\{||z|| : z\in Z_3(M;\mathcal{P},\RR),  [\psi(z)]=[M]\}$.
$||[M,\partial M]||_P$ is defined in a similar manner. 

\begin{lemma}
For $\interior M^3$ a finite volume hyperbolic manifold, 
$Vol(\interior M) = V(P) ||[M,\partial M]||_P$.
\end{lemma}
\begin{proof}
This is proven in the same manner as in \cite{BP,Gr}. 
We review the key elements. One takes a sequence of compact polyhedra $P_i$
converging to the ideal maximal $P$, and smears $P_i$ uniformly about $M$
to get a measure chain on $M$, where the copies of $P_i$ with opposite orientation
to M have negative sign. 
This measure chain is a cycle, since  any
face $F$ of $P_i$ will be matched with $\overline{F}$ as a face of the polyhedron
obtained by reflecting $P_i$ through the geodesic plane containing $F$. One
can approximate this measure cycle by choosing a Dirichlet domain and a
basepoint. Then send each vertex of a polyhedron in a measure cycle to
the nearest basepoint vertex. Each polyhedron then is weighted by the
measure of the set of polyhedra sent to it under this map. One checks
that for large enough $i$ we still have a fundamental 
$\mathcal P$ cycle of the same norm. Then
this chain has norm $Vol(M)/Vol(P_i)$, and taking the limit as $i\to\infty$
gives the desired result.

The non-compact case works similarly to Theorem 6.5.4 in Thurston's notes
\cite{Th}.
\end{proof}

Suppose that $\interior M$ is finite volume hyperbolic. 
Let $M'$ be a copy of $M$ with the opposite orientation,
and let $DM=M\cup_{\partial M=\partial M'} M'$ be the
double of $M$ along $\partial M$. Take a relative $\mathbb P$-cycle
$z$ such that $[\psi(z)]=[M,\partial M]$. Let $z'$ be
the corresponding cycle for $M'$. Since $M'$ has 
reversed orientation from $M$, $\partial z'$ has
opposite orientation from $\partial z$. So $\partial(z+z')=0$,
thus $z+z'$ is a fundamental cycle for $DM$. Thus
$||[DM]||_P \leq 2 ||[M,\partial M]||_P$.

\section{Normalization of cycles}
Let $S$ be incompressible in $M$, $P$ be a polyhedron. 
Take a $\mathcal P$-cycle $\mu =\sum_i a_i\mu_i \in H_3(M;\mathcal P,\RR)$.
A {\it normal} embedded disk in $P$ is one which meets each edge of $P$ at most
once. The cycle $\mu$ is {\it normal with respect to S} if each
$\mu_i:P\to M$ is transverse to $S$, and $\mu_i^{-1}(S)$ consists
of normal disks in $P$.

\begin{lemma} \label{normal cycles}
We can homotope $\mu$ to be a normal cycle with respect to $S$.
\end{lemma}
\begin{proof}
First, make $\mu$ transverse to S. 
Consider the universal cover
$\widetilde{M} \cong \mathbb{R}^3$. Then $\widetilde{S} \subset \widetilde{M}$ is
a disjoint union of properly embedded planes. $\mu$ has preimage a locally
finite cycle $\tilde{\mu}$ in $\tilde{M}$ which is equivariant under the
$\pi_1 M$ action. For each pair of vertices in $\tilde{\mu}$, choose
an edge connecting these points which meets $\tilde{S}$ in as few points
as possible, such that the choice is equivariant. Construct a new
cycle $\nu=\sum_i a_i\nu_i$ inductively. Map the vertices of each
$\nu_i$ to the corresponding vertices of $\mu_i$. Then extend
$\nu_i$ to $P^{(1)}$ by the choices of efficient edges.

 For a given
lift $\tilde{\nu}_i$ of $\nu_i$, and a component $\tilde{S}_0$ of
$\tilde{S}$, $\tilde{\nu}_i^{-1}(\tilde{S}_0)$ meets a cutset
of edges in $P^{(1)}$ in at most one point per edge. We  can extend these 
to normal curves in $\partial P$, and then to normal disks in $P$ 
(figure \ref{normal}). To see this, look at a face $\Delta$ of
$P$. If one edge of $\partial\Delta$ meets $\tilde{\nu}_i^{-1}(\tilde{S}_0)$,
then exactly one of the other two edges of $\partial\Delta$ meets 
$\tilde{\nu}_i^{-1}(\tilde{S}_0)$. So we may connect the two
points of $\partial\Delta\cap \tilde{\nu}_i^{-1}(\tilde{S}_0)$ 
by a normal arc in $\Delta$. We may always choose these arcs
disjoint, for different components of $\tilde{S}$. If not, then
there would be components $\tilde{S}_0$ and $\tilde{S}_1$
such that $\tilde{\nu}_i^{-1}(\tilde{S}_0)\cap\partial\Delta$
and $\tilde{\nu}_i^{-1}(\tilde{S}_1)\cap\partial\Delta$ are
linked in $\partial\Delta$. But $\tilde{S_0}$ divides
$\nu_i(\partial\Delta)$ into two pieces, which lie on different
sides of $\tilde{S_0}$. So only once piece could intersect
$\tilde{S_1}$, which means the points are not linked. 
\begin{figure} [htb]
	\begin{center}
	\psfrag{M}{\small$\tilde{M}$}
	\psfrag{S}{\small$\tilde{S}$}
	\psfrag{n}{$\tilde{\nu}_i$}
	\psfrag{P}{$P$}
	\epsfxsize=4in
	\epsfbox{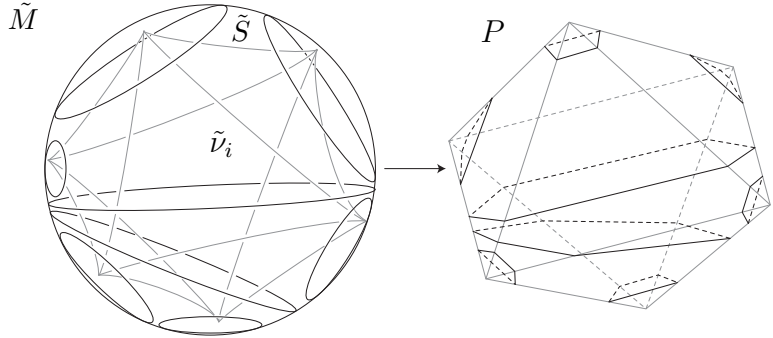}
	\caption{\label{normal} The surface S determines canonical normal disks}
	\end{center}
\end{figure}
Now, extend $\tilde{\nu_i}$ to these curves, and then to the disks
equivariantly so that the normal disks map to $\tilde{S}$. This 
can be done since $\tilde{S}$ is just a union of planes, and
the action of $\pi_1(M)$ is free. We may
now extend to $P^{(2)}$ so that the complement of the normal curves
misses $\tilde{S}$, and then we can fill in the balls in P complementary
to the normal disks equivariantly, and missing $\tilde{S}$, since
the complementary pieces of $\tilde{S}$ are contractible. This
chain $\tilde{\nu}$ is equivariantly homotopic to $\tilde{\mu}$,
so the chain $\nu$ is the desired normalized chain. 
\end{proof}

Remark: This proof should also work if we take a tight lamination \L in $M$,
in which case we can make the edges meet each plane in $\tilde{M}$ at
most once, and the rest of the argument is the same.

\section{Lower bound on volume}
In this section, we prove our main theorem \ref{main}.
This will be a consequence of the following lemma, together with a
well chosen sequence of polyhedra. 

\begin{lemma}
Let $M$ be an orientable, finite volume hyperbolic manifold. 
Let $P$ be a cell with $F$ triangular faces, which has 
a realization as a convex ideal polyhedron in \H3. Let $S\subset M$
be an incompressible surface, and let $N=M\setminus\neigh(S)$. Then 
$Vol(M)\geq -2\frac{Vol(P)}{F-2} \chi(Guts(N)).$
\end{lemma}
\begin{proof}
First, let's outline the steps involved:\newline
Step 1: Normalize the fundamental cycle.\newline
Step 2: Cut the manifold and cycle along $S$.\newline
Step 3: Lift the cut up cycle to coverings corresponding to $Guts(N)$.\newline
Step 4: Put a hyperbolic metric on $Guts(N)$ and compactify the
cusps to points.\newline
Step 5: Straighten the cycle.\newline
Step 6: Collapse to a new cycle.\newline
Step 7: Show that each triangle contributes only once
to $\partial Guts(N)$ .\newline
Step 8: Estimate the contribution of each cell in the 
original cycle to the boundary cycle.\newline
Step 9: Find a lower bound in terms of $\chi(Guts(N))$.

Step 1:
Since $M$ has finite volume, $M$ may have cusps.
In this case, replace $M$ and $S$ by their doubles along the
cusps, to get a closed manifold and surface. Denote $Guts(N)=(Q,R)$.
$Q$ has components $Q_i$  and $R_i=Q_i\cap R$. 

As an example of these definitions, consider the link shown in  
figure \ref{mutant}, which is the double of the tangle on one side
of the Conway sphere. 
Double along the link to obtain $M$. $S$ is the double of the
Conway sphere. Cutting along $S$ gives an atoroidal manifold 
which is the double of the tangle inside the sphere.
Then $N$ is two copies of the  double of the tangle, the 
characteristic submanifold  
is a regular neighborhood of the tangle, and 
$R$ is boundary of a neighborhood of the tangle. $Q$ is four copies of the tangle complement.

 Take a $\mathcal P$-cycle
$\mu = \sum_{i=1}^n a_i\mu_i$, such 
that
$[\psi(\mu)] = [M] \in H_3(M;\RR)$, and $$\|M\|_P \leq \|\mu\|_P \leq 
\|M\|_P + \epsilon.$$
By the normal cycles lemma, we may assume $\mu$ is normal. 

Step 2: For each singular cell $\mu_i : P \to M$, we can cut it
along $\mu_i^{-1}(\mathcal N(S))$, to obtain cells 
$P_i^1,...,P_i^{p_i}$, where $\mu_i : P_i^j \to N$(see Fig. 
\ref{normalchunks}).

\begin{figure}[htbp] 
	\begin{center}
	\epsfxsize=3in
	\epsfbox{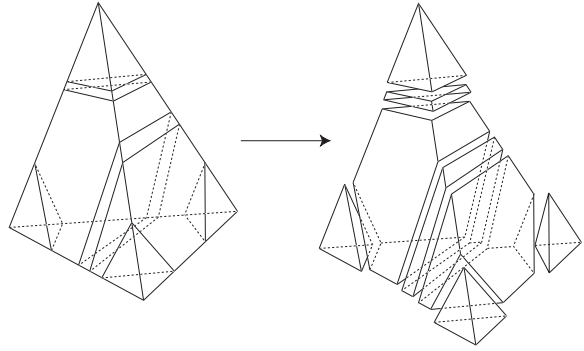}
	\caption{\label{normalchunks} Cutting up the cycle}
	\end{center}
\end{figure}

 Then we get a cell cycle 
$$\nu = \sum_{i=1}^n \sum_{j=1}^{p_i} a_i \sigma_{i|P^j_i}$$
Since [$\mu$]=[M], $\mu$ is locally degree one. So \n\ is also 
locally degree
one. We choose a decomposition of each face of 
$\partial P_i^j\cap \nu_i^{-1}(\partial N)$ into triangles, without
introducing any new vertices. 

Step 3:  Choose a component $Q_i$ of $Q$. Then $\pi_1(Q_i)$ injects into 
$\pi_1(N)$. Take
the cover $\widetilde{N}_i$ of $N$ corresponding to $\pi_1(Q_i)$.
 Then  there is a lift $Q_i \to \widetilde{N}_i$, which is 
a homotopy equivalence since it is
an isomorphism on fundamental group (we will call the image of the lift 
$Q_i$ too). There are simply
connected complementary pieces for each component of $R_i$. 
So there is a retraction $r:\widetilde{N}_i \to 
Q_i$ crushing each component of $\widetilde{N}_i\backslash Q_i$ to its 
corresponding component of $R_i$. 
The cell chain $\nu_i$ has preimage a locally finite chain
$\widetilde{\nu}_i$ in $\widetilde{N}_i$. Since $Q_i$ is compact,
we can get a finite chain $\widetilde\nu_i'$ by taking only the cells of 
$\widetilde{\nu}_i$ which intersect the lift of $Q_i$. Then we project
$\widetilde{\nu}_i'$ to $Q_i$ by the retraction $r$.

Step 4: By Thurston's geometrization theorem, $Q\setminus R$ has a 
geometric structure
of finite volume with totally geodesic boundary.  
 Crush each cusp of $R$ to a point.
This is equivalent to adding parabolic limit points to $Q$. We'll
call these parabolic points $\hat{R}$, and the new manifold
$\hat{Q}$ with cycle $\hat{\nu}$. 

Remark: The subsequent argument may be made without choosing a metric,
but it makes some later choices canonical, and gives a better
intuition for the argument, in our opinion.

Step 5: Now, we will
straighten $\hat{\nu}$ inductively. This will be done
in a certain order.  First, we straighten interior
edges. This will be done in a manner which keeps
each end point of each edge on the boundary
component on which it started (or keeps interior endpoints fixed).
If the endpoint of an edge is mapped to a point $x\in \hat{R}$,
then it will remain fixed at $x$, unless it was originally
mapped to one of the two boundary components incident
with $x$, in which case it can move around on that component.
Each edge will be homotoped to a unique arc or cusp point in this manner:
each pair of boundary components in $\widetilde{Q_i}$ has a unique
geodesic connecting them, or a unique tangent point in $\widetilde{\hat{R}}$.
By the normalization, each interior edge has end points on distinct
boundary components, when lifted to the universal cover, so
no arc will be homotopic into a component of $\partial\hat{Q}$.

Next, homotope the edges in the boundary to be
straight, keeping the endpoints fixed. Do this in such
a way that edges which homotope to the same geodesic
have the same parametrization. Finally, homotope
the triangles in each face 
in the boundary  (chosen at the end of step 2) to be straight.
We may choose these so that parametrizations of edges
are canonical, e.g. by projecting down from the
Lorentzian model of $\widetilde{\partial Q}$. 
 Call the resulting cell cycle of $\partial\hat{Q}$ 
$\alpha$. 

\begin{figure}[htbp] 
	\begin{center}
	\epsfxsize=4 in
	\epsfbox{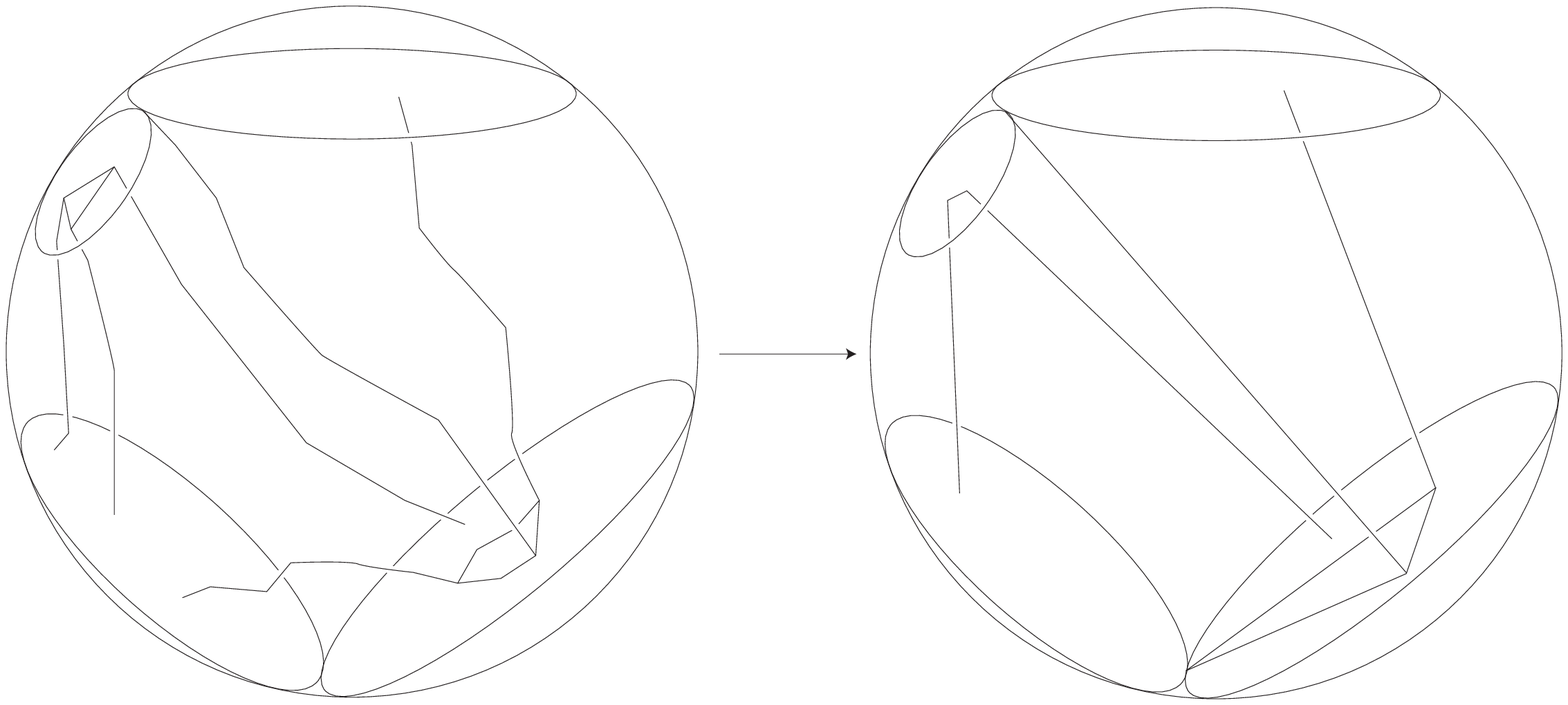}
	\caption{ Straightening the cycle}
	\end{center}
\end{figure}

We need to show that $\alpha$ is degree one on 
$\partial Q_i\cap \partial \widetilde{N_i}$. The chain $\widetilde{\nu}_i'$
is locally degree one on the interior of $Q_i$, since
$\mu$ and $\nu$ are locally degree one. Then 
$\beta=\partial\widetilde{\nu}_i\cap \partial\widetilde{N_i}$  is a
locally degree 1 chain on $\partial Q_i\cap \partial \widetilde{N_i}$.
$\partial \beta$ has support in $\partial \widetilde{N_i}\setminus
\partial Q_i$. When we retract $\widetilde{N_i}$ to $Q_i$ by
the map $r$, then all the edges of $\partial \beta$ will
lie inside $R_i$. All edges connecting points
of $\partial \beta$ will lie in $R_i$, since we chose
$\widetilde{\nu}_i$ to include every cell which intersects $Q_i$. 
So when we crush $R$ to
$\hat{R}$, and straighten, $\partial \beta$ will remain
in $\hat{R}$, so $\beta$ will remain degree one in 
$\partial \hat{Q}_i$. 

Step 6: Many triangles in $\alpha$ will have collapsed
to edges or points under this straightening process. The set of all
collapsed triangles forms a degree 0 subcycle of $\alpha$, so we may 
eliminate it and retain a degree 1 cycle in $\partial\hat{Q}$. 
We will call this new cycle $\alpha'$.

Step 7: We need to show that each triangle of $\alpha'$
contributes at most once to $\partial Q$. That is, we took the
preimage of the cycle $\nu$ in $\widetilde{N_i}$  to get $\widetilde{\nu}_i$
in step 3.
A particular triangle will have many preimages in $\widetilde{\nu}_i$,
for each $i$,
and we need to show that at most one preimage contributes non-trivially
to $\alpha'$. To see this,
first notice that the position of each triangle after straightening
(step 5) is determined by the interior edges
of the cycle $\nu$ (the original cut up cycle from step 1)
to which it was attached. Take a particular triangle $\Delta$
in $\partial N$ which has edges $e_1, e_2, e_3$ attached to
it in the particular cell of $\nu$ in which it lies. 
Now, consider the lift of $\Delta$ to $\widetilde{N}$, the
universal cover of $N$. Suppose
that $\Delta$ contributes twice to $Q$. Let $\widetilde{Q}$ be the 
preimage of $Q$ in $\widetilde{N}$. Then a lift of $\Delta$ to $\widetilde N$
is homotopic into two different components 
$\widetilde{Q}_1$ and $\widetilde{Q}_2$ of $\widetilde{Q}$. There is some 
infinite strip or plane component $\widetilde{R}_0$ of 
the preimage of $R$ (the pared locus of $Q$) in $\widetilde{Q}$ 
which separates $\widetilde{Q}_1$ and $\widetilde{Q}_2$. If 
$\widetilde{R}_0$ is a plane, then $\Delta$ can contribute 
only to the $\widetilde{Q}_i$ which lies on the same side of $\widetilde{R}_0$
as $\Delta$ does.
So $\widetilde{R}_0$ is an infinite strip, covering an annulus in $R$.
$\widetilde{R}_0$ will be incident with two components of $\partial\widetilde{N}$,
$S_1$ and $S_2$. We may assume that the endpoints of $e_1,
e_2, e_3$ all lie on distinct components of $\partial\widetilde{N}$,
otherwise $\Delta$ will always collapse. 
Let's say that $\widetilde{Q}_1$
lies to the left of $\widetilde{R}_0$, and $\widetilde{Q}_2$
lies to the right of $\widetilde{R}_0$. If $\Delta$ lies on a
component of $\partial\widetilde{N}$ to the left of $\widetilde{R}_0$,
then it will be crushed to a cusp in $Q_2$ when we 
straighten. Symmetrically for the right, so $\Delta$ must
lie on $S_1$ (or $S_2$ ) if it is to contribute to both
$\widetilde{Q}_1$ and $\widetilde{Q}_2$. We will assume $\Delta$ is on $S_1$. 
There are some cases to consider now:

\begin{itemize}

\item  $e_1$ ends on $S_2$, $e_2$ and $e_3$ end on other
components of $\partial\widetilde{N}$ to the left of $\widetilde{R}_0$.
In this case, all the edges will be homotopic to $\hat{R}$ 
in $\widetilde{Q}_2$, so $\Delta$ will not contribute to $\widetilde{Q}_2$.

\item $e_1$ ends on $S_2$, $e_2$ ends on a boundary
component on the left of $\widetilde{R}_0$  and $e_3$ ends on a boundary
component on the right of $\widetilde{R}_0$. In this case,
$e_1$ and $e_2$ will be collapsed to $\hat{R}$ in $\widetilde{Q}_2$, so 
$\Delta$ will not contribute to $\widetilde{Q}_2$ (or to $\widetilde{Q}_1$). 

\item $e_1$ and $e_2$ end on boundary components to the left
of $\widetilde{R}_0$ and $e_3$ ends on a boundary
component on the right of $\widetilde{R}_0$. In this case,
the ends of $e_1$ and $e_2$ will be sent to $\hat{R}$ in 
$\widetilde{Q}_2$, and the edge of $\Delta$ connecting $e_1$ and $e_2$ 
will collapse in $\widetilde{Q}_2$, so $\Delta$ does not contribute
(see figure \ref{triangle}).
\begin{figure}[htbp] 
	\begin{center}
	\psfrag{M}{$\widetilde{M}$}
	\psfrag{e}{$e_1$}
	\psfrag{f}{$e_2$}
	\psfrag{g}{$e_3$}
	\psfrag{S}{$S_1$}
	\psfrag{T}{$S_2$}
	\psfrag{R}{\small$\widetilde{R}_0$}
	\psfrag{D}{$\Delta$}
	\psfrag{Q}{\small$\widetilde{Q}_1$}
	\psfrag{P}{\small$\widetilde{Q}_2$}
	\epsfxsize=5 in
	\epsfbox{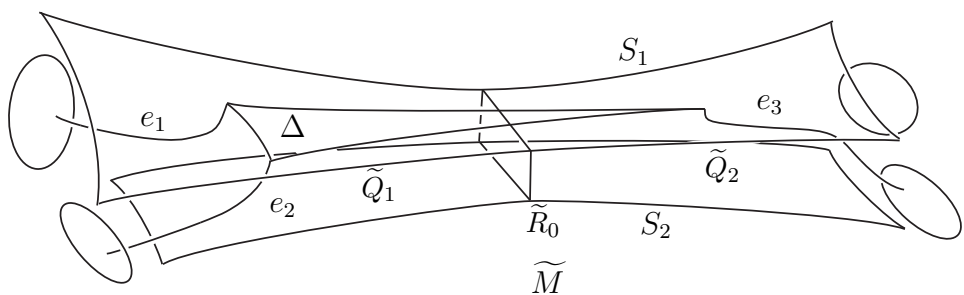}
	\caption{\label{triangle} Assessing the contribution of $\Delta$}
	\end{center}
\end{figure}
\item $e_1, e_2,$ and $e_3$ end on boundary components to the left
of $\widetilde{R}_0$.
In this case, all the ends will be mapped to $\hat{R}$ in $\widetilde{Q}_2$, 
and $\Delta$ will not contribute to $\widetilde{Q}_2$. 
\end{itemize} 

So we see that each triangle can contribute in only one
way to the cycle $\alpha'$. 

Step 8: 
Here is a way to compute the number of triangles from a
polyhedron $P_i$ which contribute to the collapsed cycle $\alpha'$.
$P_i^j$ will have faces which came from $P_i$, called $\partial P_i^{j\prime}$,
 and new faces which came
from $\partial N$, called $\partial P_i^{j\prime\prime}$. 
Collapse each quadrilateral and triangle 
in $\partial P_i^{j\prime}$ to an edge by collapsing the edges adjacent to
$\partial P_i^{j\prime\prime}$ - this collapsing is compatible with
the straightening of edges which we did in step 5 (see Figure \ref{cutting}).
 \begin{figure}[htbp] 
	\begin{center}
	\epsfxsize=3 in
	\epsfbox{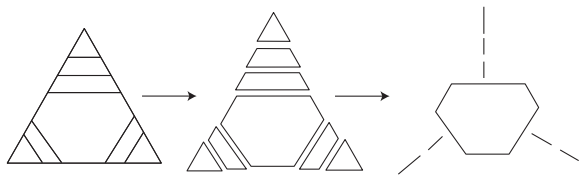}
	\caption{\label{cutting} How to collapse the faces of $P_i^j$}
	\end{center}
\end{figure}
Now, $\partial P_i^{j\prime}$ will consist of truncated triangle faces after the
collapsing (some $j$ may collapse to line segments, which we throw out). 
In step 2, each face of $\partial P_i^{j\prime\prime}$ was
divided up into triangles, and after collapsing $\partial P_i^{j\prime}$,
some of the boundary of these triangles will collapse, so we get rid
of them, since they will contribute to the collapsed cycle $\alpha-\alpha'$. 
 Now, each triangle of $\partial P_i^{j\prime\prime}$
which contributes to the cycle $\alpha'$ (from step 6) must come from
one of these triangles left over. So we need to estimate
how many of these triangles there are. Collapse the
faces of $\partial P_i^{j\prime}$ to points.
 We then get a set of triangulations of boundaries of balls, one
for each $j$ (except for the segments which we got rid of).
 Each new vertex corresponds to a cell which is divided up
into triangles of the collapsed
cycle $\alpha'$ constructed in step 6. 
If the vertex has degree $d$, then it contributes at most $d-2$
triangles to $\alpha'$. The union of triangulations has the same
number of faces as the original triangulation of $\partial P_i$, but there
are more components than the original. We may as well assume
every vertex contributes, since this is the maximal possible
case. Let there be $v$ vertices, where vertex $i$ has degree $d_i$,
$e$ edges, and $c$ components to the triangulation. 
Then $\sum_i (d_i-2) = 2 e -2 v = 2 F - 4 c\leq 2F-4 $, by euler characteristic.
So there are at most $2F-4$ triangles contributing to 
$\alpha'$ from the polyhedron $P$.

Step 9: The cycle $\alpha'$ is locally degree one. In the
metric on $\partial Q$, each triangle has area at most 
$\pi$. We have 
$$Area(\partial Q)/\pi=-2\chi(\partial Q)=-4\chi(Q)
\leq ||\alpha'|| \leq (2F-4) ||\mu|| \leq (2F-4)(||M||_P+\epsilon).$$
Letting $\epsilon$ go to zero, then 
$$Vol(M)=||M||_P Vol(P)\geq \frac{-4\chi(Q)}{2F-4}Vol(P)=
\frac{-2\chi(Q)}{F-2}Vol(P).$$
\end{proof}

To finish the main theorem, we need to find a sequence of
polyhedra $P_n$ such that $\frac{Vol(P_n)}{F_n-2}$ approaches $V_3$.
First note that this is the best one can do. An ideal polyhedron $P$ with 
$F$ triangle faces has $Vol(P)\leq V_3(F-3)$, since we can triangulate
$P$ by coning off to a vertex, and we need at most $F-3$ 
tetrahedra, each of volume $\leq V_3$.
 Sleator, Tarjan, and Thurston showed that one could find a 
sequence of polyhedra $P_n$ 
such that $\frac{Vol(P_n)}{F_n-2}\to V_3$ by taking polyhedra
which have the combinatorics of a hexagonal subdivision
of an icosahedron (this was pointed out to me by Thurston) \cite{STT}.
The idea is that most of the triangles will look
nearly equilateral when viewed from one ideal vertex, 
so that coning off will give a triangulation where
most of the tetrahedra are nearly regular. 

\begin{figure}[htbp] 
	\begin{center}
	\epsfxsize=2 in
	\epsfbox{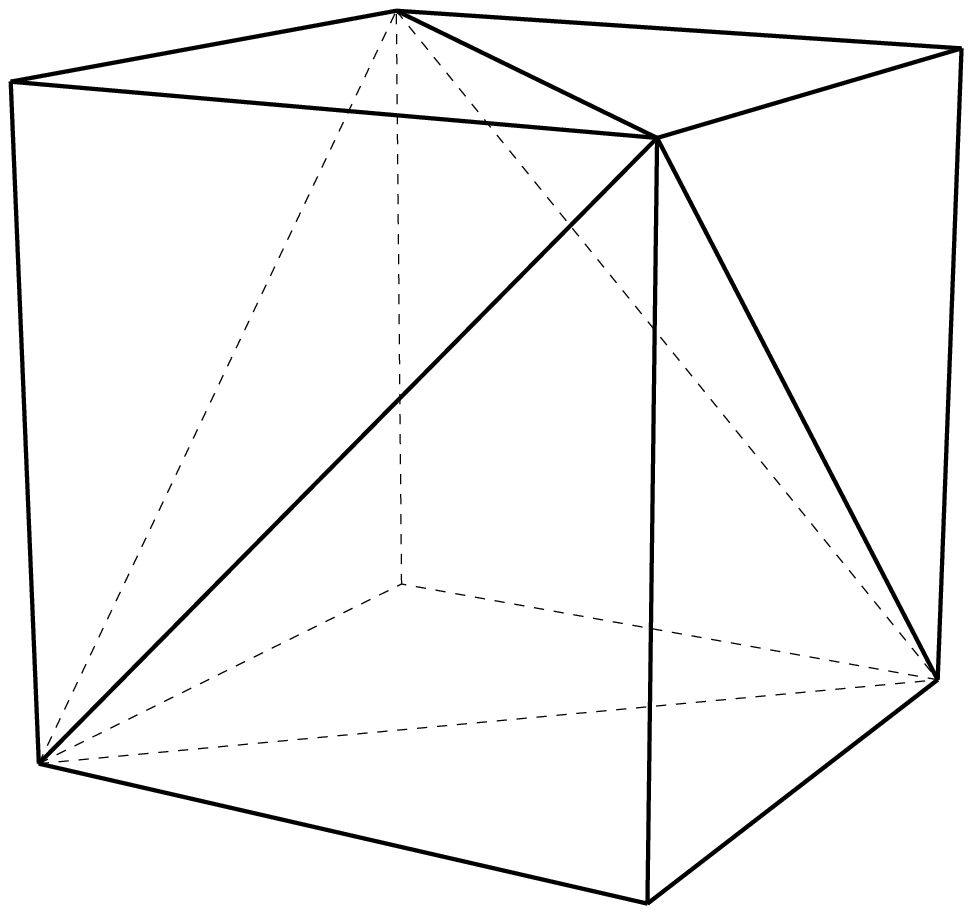}
	\epsfxsize=2.5in
	\epsfbox{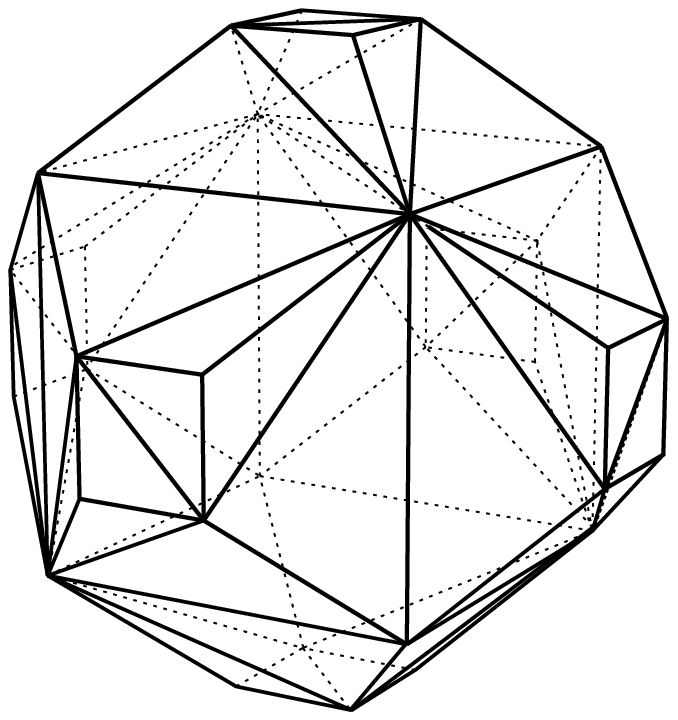}
	\caption{ \label{polyhedra} Projective model of $P_1$, $P_2$, and $P_3$}
	\end{center}
\end{figure}

Another possible choice of polyhedra is obtained by a sequence
of polyhedra which exhaust the tesselation of  \H3  by 
regular ideal tetrahedra. $P_1$ is a tetrahedron. $P_n$ has
all dihedral angles either $\pi$ or $\frac{\pi}{3}$, so it
will have triangular faces, but geometrically there will be
larger faces consisting of coplanar triangles. $P_{n+1}$ is obtained by 
adding a reflected copy of $P_n$ to each geometric 
face of $P_n$ (see figure \ref{polyhedra}). One can compute that in the limit 
the $Vol(P_n)/F_n \to V_3$. This is a more elementary
construction, but it may be less intuitive. 

\section{Two-bridge links} \label{sec2bridge}
In this section, we prove the claim in example \ref{twobridge}.
First, we need a general lemma, which allows one to apply the
main theorem. 

A graph is $n$-connected, if removing $n$ vertices from the
graph, as well as their adjacent edges, keeps the graph 
connected.

\begin{lemma}\label{handle}
Let $H$ be a handlebody of genus $>1$, $P\subset\partial H$ a set
of essential simple closed curves, and $D$ a set of properly embedded 
disks in $H$ which cut $H$ up into a ball. Take $D$ transverse
to $P$ with minimal intersection. Let $S=\partial(H\setminus\neigh(D))$
be the resulting sphere. Take a graph $G$, with
vertices coming from $\partial\overline{\neigh(D)}\cap S$ and
edges the segments of $P\cap S$. If $G$ is 2-connected, then
$\chi(Guts(H,P))=\chi(H)$. 
\end{lemma}

\begin{proof}
Since $H$ is a handlebody, it is atoroidal. So the only 
Seifert pieces in the characteristic sub-pair of $(H,P)$ must
be solid tori. Suppose the characteristic $I$-bundle
of $(H,P)$ has $\chi<0$. Then there is a sub-bundle which
is $T=t\times I$, where $t$ is a pair of pants. 

First, notice that the disks $D$ must be $\partial$-incompressible,
that is there is no $\partial$-compression with one boundary arc
on $D$ and the other boundary arc on $S\setminus G$. 
If not, then the 
$\partial$-compression would represent a separating vertex of $G$. But the
graph $G$ is 1-connected, since it is 2-connected, and has at
least 4 vertices, a contradiction. This is a diskbusting 
criterion of Whitehead. 

Look at the intersection of $D$ with the product pair of
pants $T$. Make $T\cap D$ minimal. We need to show that 
$T\cap D$ consists of product rectangles. The intersection will have no closed
curves, because $\partial t\times I$ is incompressible,
and $D\cap (\partial t\times I)$ will consist of 
product lines, by the $\partial$-incompressibility of $D$. 
Suppose that some component $A$ of $T\cap D$ is not a product rectangle
in $T$. Then there are two arcs of $\partial A$ running over $t\times \{0\}$.
Taking a path connecting these arcs on $A$, we see that this
path is homotopic rel endpoints into $t\times \{0\}$, so $D$
has a boundary compression by a well known argument, which we 
outline. Make the homotopy boundary compression transverse to $D$, with no
fake branching. Then an innermost arc of intersection on the
boundary compression gives a disk with interior disjoint from $D$. 
The loop theorem allows us to replace it with an embedded 
$\partial$-compressing disk. So there is a $\partial$-compression
of $D$ inside $T$. But since $D$ is $\partial$-incompressible 
in the complement of $G$, when we surger along the $\partial$-compression,
we get something isotopic to $D$, with smaller intersection with
$\partial t\times I$, a contradiction. So $D$ must be a product
rectangle.

If we cut $T$ along $T\cap D$, we must get a set of disks $\times I$,
otherwise there would be a simple closed curve in the sphere
which compresses $t\times\partial I$, a contradiction. Call one
of these disks $d\times I$. $\partial d\times I$ intersects $D$
in a collection of at least three rectangles. Each of these
rectangles intersects two vertices of $G$ and cuts the graph
into two pieces. Since $G$ is 2-connected, the part of $G$ 
lying to one side of each rectangle must contain no vertices of $G$. 
So we must see a picture like figure \ref{product}. 
\begin{figure}[htb]
	\begin{center}
	\psfrag{G}{$G$}
	\psfrag{P}{$P$}
	\psfrag{T}{$d\times {0}$}
	\psfrag{D}{$D$}
	\psfrag{R}{$d\times {1}$}
	\epsfxsize=1.5in
	\epsfbox{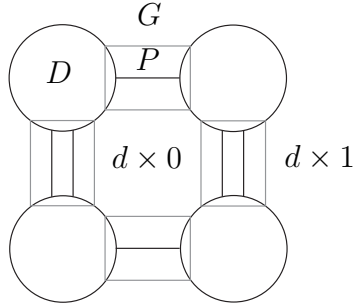}
	\caption{\label{product} How a product region must look}
	\end{center}
\end{figure}
Thus,
$d\times I$ must meet every vertex in $G$. This picture is 
impossible, since it is not 2-connected. 
\end{proof}

We need another observation, which is essentially due to Gabai
in the context of sutured manifolds \cite{Ga2}. If in the context of the
lemma, a disk hits the curve exactly twice, then this gives
a pair of pants $\times I$. This means that if we have 
the graph $G$ as in lemma \ref{handle}, we can get rid
of any degree 2 vertices, since they contribute to the
characteristic submanifold. If what's left over is 
2-connected, then we may apply lemma \ref{handle} to each
component which is left over. 

We can now apply the lemma to analyze volumes of 2-bridge links.
We will follow the conventions of Hatcher and Thurston \cite{HT}.
A $2$-bridge link $L$ is obtained by taking the boundary of
a plumbing $S$ of several bands $B_i$ (twisted moebius strips and
annuli), each having $b_i$  half-twists, where $|b_i|\geq 2$,
and the twists are right-handed if $b_i>0$, and left-handed
if $b_i<0$.

 To estimate the volume, consider the sequence of
numbers of twists in each band $\{b_1,\dots,b_k\}$. Consider
the number $j$ of maximal subsequences $b_{k},b_{k+1},...,b_{k+l}$, where
each $|b_i|\geq 3$, and let $m = |\{i : |b_i|\geq 3\}$. 
Then $Vol(L)\geq 2 V_3 (m-j)$. Each sequence of bands counted 
by $j$ contributes to the guts of the complement of the surface.
Thicken the surface $S$ up to get a Heegaard handlebody,
with $L$ on its boundary. Then choose the obvious set of
disks $D_i$, one for each band $B_i$, getting a graph $G$ as in the
lemma, with two vertices for each disk of $D$, and edges
coming from the pieces of $L$. First, let's see what the
graph $G$ looks like locally for each band. 

\begin{figure} [htb] 
	\begin{center}
	\epsfxsize=.4in
	\epsfbox{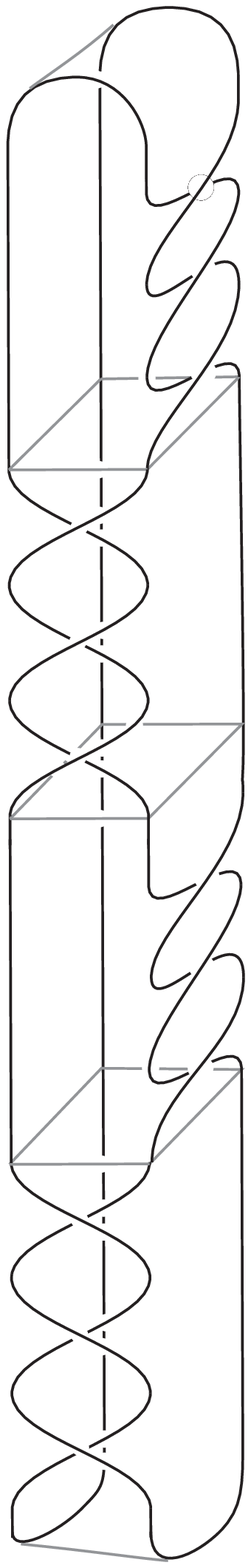}
	\psfrag{D}{\footnotesize D} 
	\psfrag{B}{\tiny $B_i$}
	\psfrag{a}{(a)}
	\psfrag{b}{(b)}
	\psfrag{c}{(c)}
	\psfrag{d}{(d)}
	\psfrag{m}{\scriptsize $b_i$}
	\psfrag{t}{\scriptsize twists}
	\psfrag{q}{\scriptsize twist}
	\psfrag{n}{\tiny $b_i-2$}
	\epsfxsize=4in
	\epsfbox{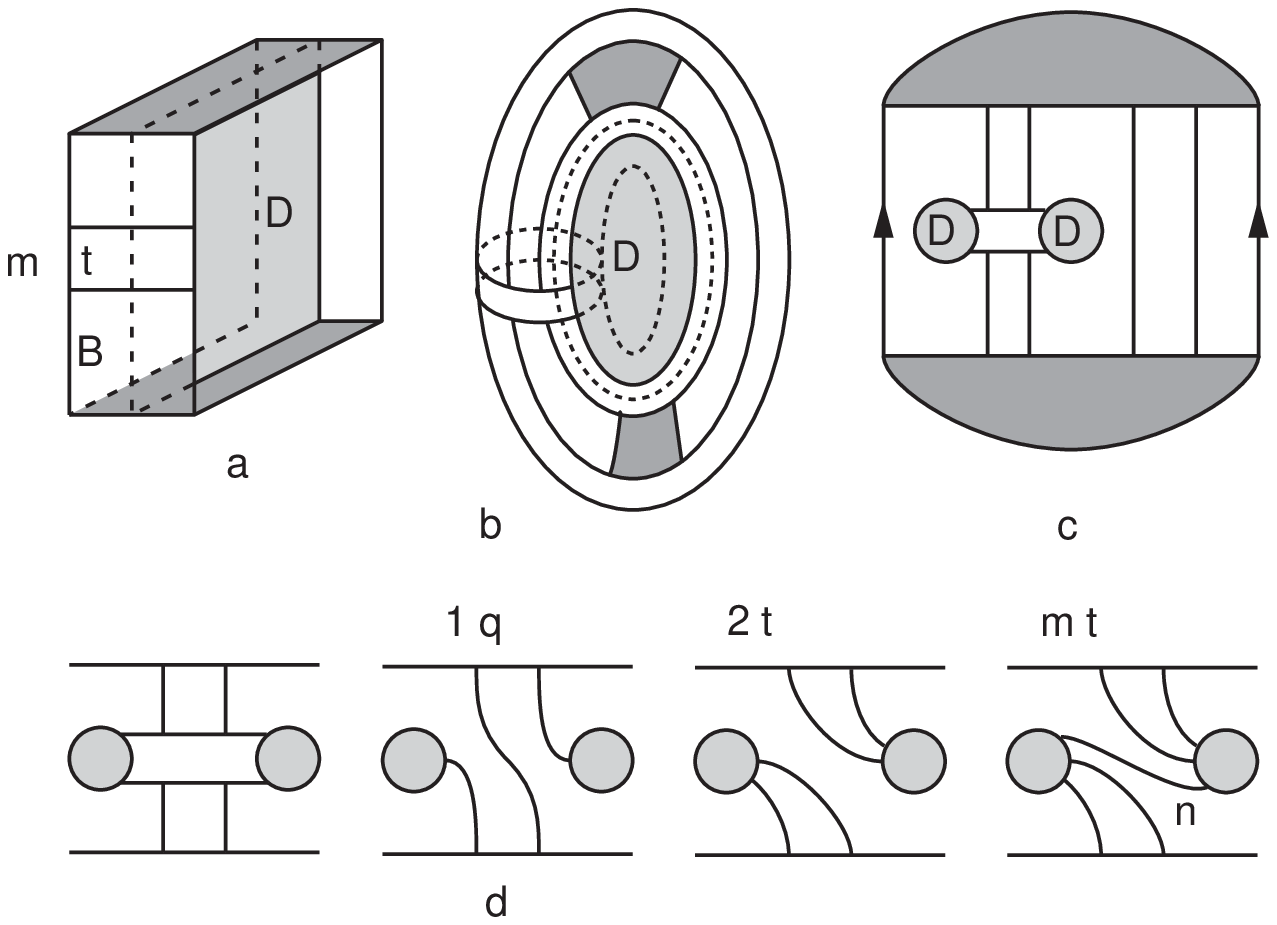}
	\caption{\label{2bridg} Constructing $G$}
	\end{center}
\end{figure}

Figure \ref{2bridg}(a) shows what one of the bands looks like,
with rectangular pieces on top and bottom where the plumbing
occurs. Fattening up and adding the 2-handle corresponding
to $D_i$, we get a sphere as in figure \ref{2bridg}(b). 
The graph on the surface of the sphere looks like figure \ref{2bridg}(c),
where the rectangle has a picture looking like one of the 
cases in figure \ref{2bridg}(d).

\begin{figure} [htb] 
	\begin{center}
	\psfrag{a}{\tiny $b_1-2$}
	\psfrag{b}{\tiny $b_2-2$}
	\psfrag{c}{\tiny $b_3-2$}
	\psfrag{d}{\tiny $b_k-2$}
	\psfrag{e}{\dots}
	\epsfxsize=\textwidth
	\epsfbox{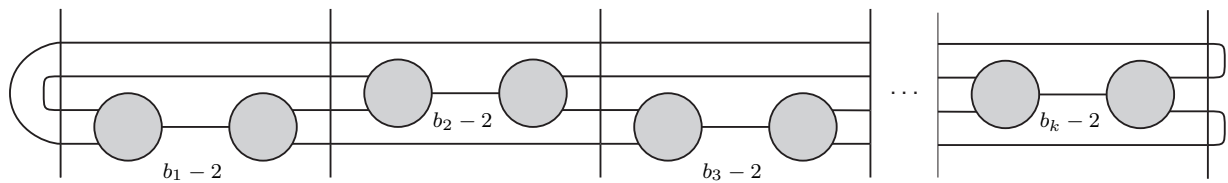}
	\caption{\label{graph2} The full graph}
	\end{center}
\end{figure}

 When we plumb several bands together, 
we get a graph which looks like figure \ref{graph2}

For each band $B_i$ with $|b_i|=2$, we can add another handle
which crosses $L$ twice, dividing the sphere into two spheres,
and creating a graph $G'$. 
Since both $D_i$ and the second handle cross $L$ twice, 
we may eliminate them from the graph $G'$. Doing this for
every band with $|b_i|=2$ leaves over a bunch of handlebodies
corresponding to maximal subsequences of $\{b_i\}$ with
$|b_i|\geq 3$. Then we may apply the lemma \ref{handle}
to these handlebodies, to obtain the claimed lower bound for
the volume. 

Let's see what this implies for 2-bridge links in general.
Any 2-bridge link may be represented by a rational number
$\frac{p}{q}$, with $0 < \frac{p}{q} < 1$. Then $\frac{p}{q}$ has
a unique continued fraction expansion

\begin{equation}
\frac{p}{q}=\cfrac{1}{a_i+
		\cfrac{1}{a_2+
		\cfrac{1}{\ddots +
		\cfrac{1}{a_k}}}}
\end{equation}
where $a_i\geq 1$, $a_1\geq 2$, $a_k\geq 2$. 

Corresponding to this expansion is a picture
in the arc complex of the 4-punctured sphere, which looks
like:

\begin{figure} [htb] 
	\begin{center}
	\psfrag{a}{\dots}
	\psfrag{b}{\footnotesize $a_1$}
	\psfrag{c}{\footnotesize $a_2$}
	\psfrag{d}{\footnotesize $a_3$}
	\psfrag{e}{\footnotesize $a_4$}
	\psfrag{f}{\footnotesize $a_{k-1}$}
	\psfrag{g}{\footnotesize $a_{k}$}
	\psfrag{h}{\dots}
	\psfrag{p}{$\frac{p}{q}$}
	\psfrag{i}{\footnotesize $\infty$}
	\psfrag{1}{\footnotesize $\frac01$}
	\psfrag{q}{$\frac{p}{q}$}
	\epsfxsize=4in
	\epsfbox{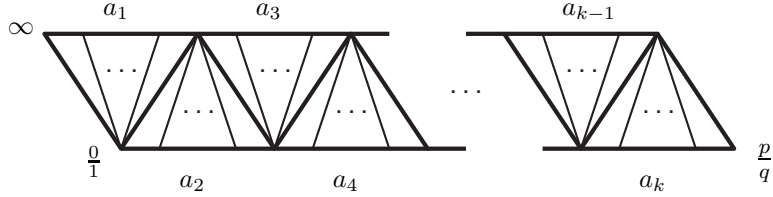}
	\caption{\label{cfrac} Continued fraction expansion}
	\end{center}
\end{figure}

A locally minimal path gives a continued fraction expansion
for $\frac{p}{q}$ which corresponds to a plumbing of bands
with $|b_i|\geq 2$, where at the $i$-th vertex, the
path turns left or right across $|b_i|$ triangles,
left if $b_i>0$, right if $b_i<0$. This must be a path traversing only heavy
edges in figure \ref{cfrac}. One can check that there is always
some path with a sequence of at least two $B_i$'s with 
$|b_i|\geq 3$, except in two classes of exceptional 
cases.
\begin{figure} [htb] 
	\begin{center}
	\psfrag{a}{\dots}
	\psfrag{b}{}
	\psfrag{c}{\footnotesize $a_2$}
	\psfrag{i}{\footnotesize $\infty$}
	\psfrag{1}{\footnotesize $\frac01$}
	\psfrag{q}{$\frac{p}{q}$}
	\epsfxsize=3in
	\epsfbox{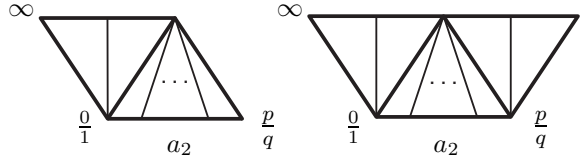}
	\caption{Cases where lower bound fails}
	\end{center}
\end{figure}
The first case corresponds to the twist knots, which are 
fillings on the Whitehead link and its mutant, so they have
volume $<V_{oct}$. The second case corresponds to two-component
links which are branch covers over twist knots in $S^3$ or 2-bridge
knots in $\RR\PP^3$. 

For every other 2-bridge link, we get that the volume $\geq 2V_3$.
Of course, this is already implied by Cao and Meyerhoff's result
\cite{CM}, but if conjecture \ref{conj1} is true, then
we could improve this lower bound to $V_{oct}$.
So conjecture \ref{conj1} would imply
that the only 2-bridge links with volume $\leq V_{oct}$ are
twist knots.

\section{Meridianal planar incompressible surfaces}\label{planar}
Let $k\subset M$ be a hyperbolic link, $M_k=M\setminus\neigh (k)$,
and let $P\subset M_k$ be an incompressible planar surface
with boundary components being meridians of $k$. $P$ extends
naturally to a $2$-sphere $\hat{P}$ in $M$, by capping off $\partial P$
with meridian disks in $\neigh (k)$. If $M$ is a $\ZZ_2-$homology
$3$-sphere, then $M_{\hat{P}}=M\setminus\neigh(\hat{P})$ is
two $\ZZ_2-$homology $3$-balls, and $k\cap M_{\hat{P}}$ is
an atoroidal tangle in each component of $M_{\hat{P}}$. For
example, a Conway sphere in a link complement is of this
form.

\begin{lemma} Let $B$ be a $\ZZ_2-$homology ball, $T\subset B$
a tangle in $B$. Denote $B_T=B\setminus\neigh(T)$, and 
$R=\partial B_T\setminus \interior\partial B$. Assume that $\partial B\cap 
\partial B_T$ is incompressible in $B_T$
and $B_T$ is atoroidal. Then
$\chi(Guts(B_T,R))\leq -1$. 
\end{lemma}
\begin{proof} Let $(G,\partial_0 G)=Guts(B_T,R)$.
If $\chi(G)=0$, then $G$ is a union of $S^1\times I\times I$. 
The $I$-bundle part of the characteristic
sub-pair of $(B_T,\partial B_T\setminus \interior R)$
must have planar base surface, since $\partial B$ is a 2-sphere.
  Also,
it must be a product bundle, since otherwise there would be
a properly embedded  m$\ddot{\text{o}}$bius band in $B$, contradicting that
$B$ is a $\ZZ_2$-homology ball.
The other components of the characteristic sub-pair are solid
tori, such that the meridian meets $\partial B$ at least
3 times. We will denote these pieces by $(C,\partial_0 C)$, where
$\partial_0 C =\overline{\partial C\setminus\partial B}$. 
Fill in the components of $G$ to get a book of $I$-bundles.
Notice that $\partial B_T\setminus R$ is a 
connected surface, since its complement in $\partial B$ is a union
of disks. 

Take a path $q$ in $\partial B_T\setminus\partial_0 B$ which connects points
$q_1$ and $q_2$ on opposite sides of a product piece of 
$B_T\setminus C$, with the property that
this path hits the annuli in 
$\partial C\cap \partial B$ a minimal
number of times. This number must be $\geq 1$, since the piece is
a product. Consider the first annular    piece of 
$\partial C\cap \partial B$ which the path
crosses. It will first cross the boundary of an annulus 
$A_1\subset \partial_0 C$, then an annulus $A_2\subset \partial_0 C$,
and then  proceeds back into the product part. The path must hit these 
components of $\partial A_1$ and $\partial A_2$ at most once, 
otherwise we could find
a path hitting $\partial C$ fewer times. Now, suppose the path
crosses the other component of $\partial A_2$ at some point. Then
we could find a subpath which connects opposite points of $\partial A_2$ 
and which intersects $\partial_0 C$ fewer times (see \ref{graph} (a)). 
\begin{figure} [htb]
	\begin{center}
	\psfrag{B}{$A_2$}
	\psfrag{A}{$A_1$}
	\psfrag{r}{$q_1$}
	\psfrag{t}{$q_2$}
	\psfrag{G}{$C$}
	\psfrag{q}{$q$}
	\psfrag{u}{$q'$}
	\psfrag{a}{(a)}
	\psfrag{b}{(b)}
	\epsfxsize=\textwidth
	\epsfbox{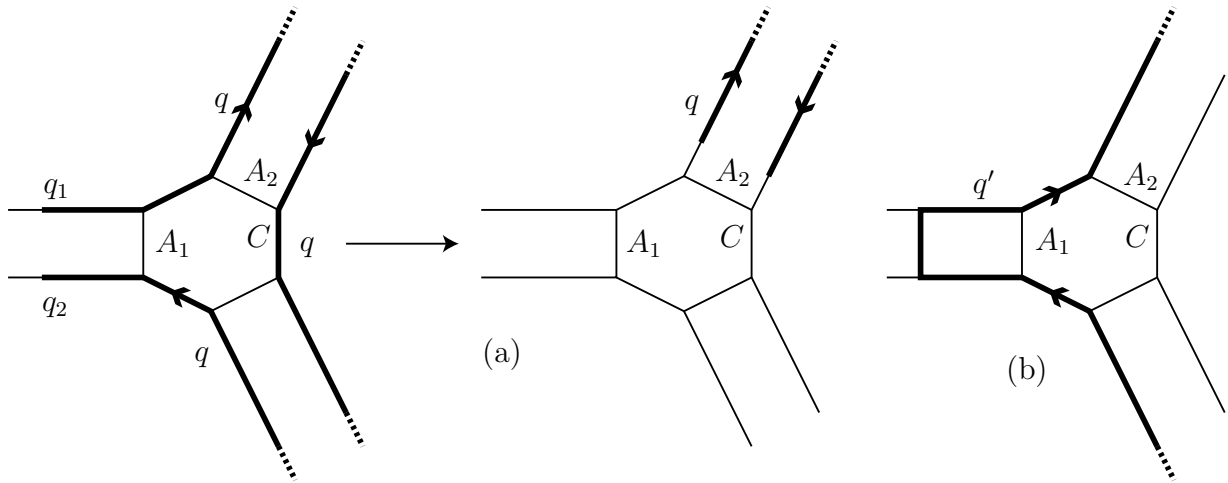}
	\caption{\label{graph} modifying the path}
	\end{center}
\end{figure}
So we may assume
that $q$ hits $\partial A_2$ at most once. Make a closed path $q'$ by 
connecting $q_1$ and $q_2$ by a path in the product part (see \ref{graph} (b)). Then
we have a closed path in $B$ which hits the annulus $A_2$ 
once, contradicting that $B$ is a $\ZZ_2-$homology ball.
Thus, $\chi(G)\leq -1$. 
\end{proof}

\begin{cor}
Let $k$ be a link in $M$, a $\ZZ_2$-homology 3-sphere with
a meridianal incompressible surface $P$. Then $Vol(M_k)\geq 4V_3$.
\end{cor}
\begin{proof}
By the lemma, $\chi(Guts(M_k\setminus\neigh(P),\partial M_k\setminus
\neigh(P)))\leq -2$, so $Vol(M_k)\geq 4V_3$.
\end{proof}

A theorem of Thompson shows that every 
knot in $S^3$ either has a meridianal incompressible planar
surface, or the knot is in bridge position whenever it is
in thin position \cite{Tho}. So this shows that many knots have thin
position equal to bridge postion, since their volumes are
less than $4 V_3$. There are 23 knots in the census \cite{CDW} which have
volume $\leq 4V_3$.

\section{Betti number 1 case} \label{cyclic2}
\begin{figure}[htbp] 
	\begin{center}
	\psfrag{Z}{$\widetilde{M}$}
	\psfrag{a}{(a)}
	\psfrag{b}{(b)}
	\psfrag{A}{$S$}
	\psfrag{S}{$S_0$}
	\psfrag{T}{$S_1$}
	\psfrag{R}{$M_1$}
	\psfrag{L}{$M_j$}
	\psfrag{U}{$S_2$}
	\psfrag{V}{$S_j$}
	\psfrag{W}{$S_{j+1}$}
	\psfrag{X}{$S_{j+1}$}
	\psfrag{t}{$t$}
	\psfrag{P}{\small$\Phi_{j+1-}^p$}
	\psfrag{K}{$\Phi_{j-}^p$}
	\psfrag{B}{$t(\Phi_{j-}^p)$}
	\psfrag{J}{}
	\psfrag{M}{$M$}
	\epsfxsize=5 in
	\epsfbox{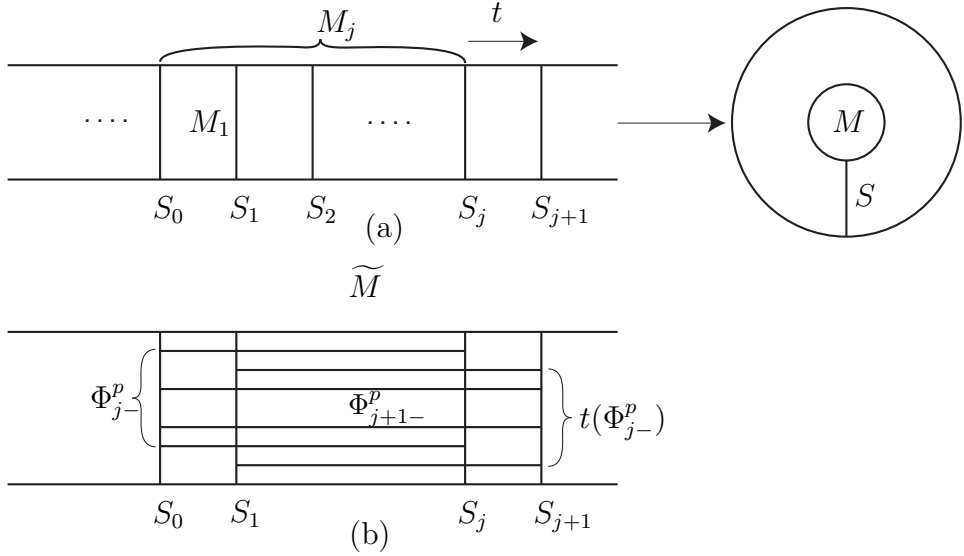}
	\caption{\label{cyclic} Schematic for the proof}
	\end{center}
\end{figure}

Consider a closed hyperbolic manifold $M^3$ with a non-separating
incompressible 
surface $S$, which does not fiber over the circle. $S$ represents
a non-trivial class in $H_2(M;\RR)$. As noted before,
a theorem of Thurston implies such a surface exists if 
$\beta_1(M)\geq 2$. We
modify an argument of Cooper and Long \cite{CL} to get a lower bound
for the volume of $M$. 
Let $\tilde{M}_n$ be the $n$-fold cyclic cover of $M$ dual to $S$, and 
$\tilde{M}$ be the infinite cyclic cover. Choose a base 
surface $S_0$ in $\tilde{M}$, let $t$ be the generating covering 
translation in $\tilde{M}$.
Let $S_i=t^i(S_0).$ In $\tilde{M}$, let $M_j$ be the manifold
between $S_0$ and $S_j$, and let $(\Sigma_j,K_j)\subset (M_j,\partial M_j)$
be the characteristic subpair, $\Phi_j\subset \Sigma_j$ be the characteristic
$I$-bundle. $M_j$ is atoroidal, so $\Sigma_j\setminus \Phi_j$ is a 
disjoint union of solid tori, which intersect $\partial M_j$ in 
a collection of essential annuli. 
$\Phi_j$ has components $\Phi_j^0$ 
whose boundary is entirely in $S_0$, $\Phi_j^j$ whose 
boundary is in $S_j$, and 
$\Phi_j^p$ which is a product between $S_0$ and $S_j$ 
(see figure \ref{cyclic}(a)).
\begin{lemma}
If $S$ is Thurston norm minimizing, then  
$\chi(Guts(M_{5g-5})) = \chi(S)=2-2g$.
\end{lemma}

\begin{proof}
Put a partial order on subsurfaces of $S_0$, such that if $F, G$
are subsurfaces of $S_0$, then $F\leq G$ if and only if $F$ 
is isotopic to a subsurface of $G$. Let $\Phi^p_{j-}$ be
the union of components of $\Phi_j^p$ with $\chi <0$.  
Let $\phi_j=\partial\Phi_{j-}^p\cap S_0$.
For a collection of curves $C$ in $S_0$, let $[C]$ represent the
set of homotopy classes of the curves. 

\begin{claim}
$\chi(\Phi_j^0)=\chi(\Phi_j^j)=0$. 
\end{claim}

\begin{proof}
By a result of Gabai, $S_0$ is Thurston-norm minimizing
in $\tilde{M}_j$ \cite{Ga3}. 
If $\chi(\Phi_j^0)<0$, then $\partial \Phi_j^0 =F\cup J$, where
$F= S_0\cap \Phi_j^0$, and $J$ is a union of annuli. Then 
$(S_0\setminus F)\cup J$ is a surface homologous to $S_0$ with
smaller Thurston norm, a contradiction. Similarly, $\chi(\Phi_j^j)=0$.
(Note: Gabai's result simplifies
this argument, but is not mandatory.)
\end{proof}

\begin{claim}
$\phi_j$ is a decreasing sequence, and
for each $j$, either $[\partial\phi_j] \subsetneq [\partial\phi_{j+1}]$,
or $\chi(\phi_j) > \chi(\phi_{j+1})$.
\end{claim}

\begin{proof}
$\Phi_{j+1-}^p$ can be isotoped so that $S_j\cap \Phi_{j+1-}^p$
is incompressible, and divides each component of $\Phi_{j+1-}^p$ into product
pieces (see Jaco, IX.1, \cite{Jac}).
The union of first pieces of each component is a product from 
$S_0$ to $S_j$.
By the characteristic submanifold theory, this product part
is isotopic into $\Sigma_{j}$, so since each component has
$\chi <0$, it must be isotopic into $\Phi_{j-}^p$ (they
couldn't be isotopic into a solid torus component of $\Sigma_j$) (fig. \ref{cyclic}(b)). Thus
$\phi_{j+1}\leq \phi_j$. By an isotopy, we may assume that
$\phi_{j+1}\subset \interior \phi_j$.

Suppose that $\chi(\phi_{j+1})=
\chi(\phi_j)$. Then $\phi_j\setminus\interior\phi_{j+1}$ is
a disjoint union of annuli. For each curve in $\partial\phi_j$,
the other boundary curve in the annulus in  
$\phi_j\setminus\interior\phi_{j+1}$ must be a curve in
$\partial\phi_{j+1}$, otherwise $\phi_j$ would have an annulus
component, contradicting that each component has $\chi<0$. 
So we have $[\partial\phi_j]\subset[\partial \phi_{j+1}]$.

Now, assume that $\chi(\phi_{j+1})= \chi(\phi_j)$ and
$[\partial\phi_j]=[\partial \phi_{j+1}]$. Then 
$\phi_{j+1}$ is isotopic to $\phi_j$.
We can isotope $\Phi_{j+1-}^p$ so that $S_1$ dissects
$\Phi_{j+1-}^p$ into product pieces. Let $P$ be the union
of these pieces which have a boundary component on $S_{j+1}$. Then $P$ is
a product from $S_1$ to $S_{j+1}$. $P$ can be isotoped into 
$t(\Phi_{j-}^p)$ (fig. \ref{cyclic}(b)). So $F=t^{-1}(\partial P\cap S_1)$ will be 
a subsurface of $\phi_j$. $F$ is homeomorphic to $\phi_{j+1}$. 
Then as with $\phi_{j+1}$, $F$ is isotopic to $\phi_{j}$.
Thus, $[\partial F]=[\partial \phi_j]$. 

Since $[\partial \phi_j]=[\partial \phi_{j+1}]$, we may take
each component of $\partial\phi_j$ to a parallel component of
$\partial\phi_{j+1}$. Then $\overline{\Phi_{j+1-}^p\setminus P}$
gives a product from $\phi_{j+1}$ to $t(F)$. So we may map
each component of $\partial\phi_{j+1}$ to a component of $\partial t(F)$
which cobounds an annulus with it, coming from the boundary
of the product. Then this component of 
$\partial t(F)$ corresponds to a component of $\partial\phi_j$, since
$\partial F=\partial\phi_j$. So we have a map from 
$[\partial \phi_j]\to[\partial\phi_j]$. Some iterate of this
map has  a fixed curve. The sequence of annuli connecting the
iterates of this curve give a torus in $M$ which is non-trivial,
contradicting that $M$ is hyperbolic. Thus we see that
if $\chi(\phi_j)=\chi(\phi_{j+1})$, then 
$[\partial\phi_{j+1}]\subsetneq[\partial\phi_j]$, which
is equivalent to the last part of the claim.
\end{proof}

\begin{claim}
The sequence $\{\phi_j\}$ can have length at most $5g-5$.
\end{claim}
\begin{proof}
We can assume that we have nested
sequence of surfaces $\phi_j$ in $S_0$, 
such that for each $j$, either $[\partial\phi_j]\subsetneq[\partial\phi_{j+1}]$,
$\chi(\phi_{j+1})>\chi(\phi_j)$. Let $C=\cup_j [\partial\phi_j]$. Then
$|C|\leq 3g-3$. Each curve of $C$ can be in $[\partial\phi_{j+1}]\setminus
[\partial\phi_j]$ at most once, since once a curve of $C$ disappears
from $[\partial\phi_j]$, it never appears later in $[\partial\phi_k]$,
$k>j$. Also, $\chi(\phi_j)$ can increase at most $2g-2$ times. 
Thus, $\phi_{5g-5}=\emptyset$. 
\end{proof} 

So $\chi(Guts(M_{5g-5}))=\chi(M_{5g-5}\setminus \Sigma_{5g-5})
=\chi(M_{5g-5}\setminus \Phi^p_{5g-5-})=\chi(M_{5g-5})=\chi(S)$.
\end{proof}

Now, we have that $(5g-5)Vol(M)=Vol(\tilde{M}_{5g-5})\geq 2V_3(2g-2)$.
Thus, $Vol(M)\geq \frac45 V_3$. So we get a universal lower bound
to volumes of these special types of manifolds.

We believe that this lower bound can be improved in general. 
In particular, if $g=2$, then we can show that $Vol(M)\geq V_3$,
since the covering argument can be improved from 5 to 4 by
a case by case analysis.

\section{Conclusion}

There are many questions which have arisen
during this paper. One is whether this
argument can be extended to genuine laminations.
This seems possible, at least if the lamination
is tight. Other possible extensions would be
to apply the technique to get lower bounds
for two-bridge links to other classes of
manifolds where the incompressible surfaces are
well understood, {\it e.g.} torus bundles over $S^1$
and alternating links. It seems likely that
arguments in section \ref{cyclic2} are not 
sharp, so it would be interesting to improve
on them. Of course, the main question is to
try to prove conjecture \ref{conj1}, which would
give several sharp lower bounds, as explained. 
The technique in this paper probably can't
be generalized to prove the conjecture, if it
is true, although one might be able to
show that one can replace $2V_3$ with $V_{oct}$
in theorem \ref{main}
using suitable generalizations of this technique. 
The only evidence we have for conjecture \ref{conj1}
is from playing with SnapPea \cite{W}, where
if you replace the link shown in figure \ref{mutant}
with different braids connecting up the tangles
on either side of the sphere, the volume seems to
go up. We can also show that $Vol(M)\geq \frac{V_3}{V_{oct}}
Vol(Guts(M\setminus \neigh(S)))$. This follows
using the technique of theorem \ref{main}, since
the volume of a truncated tetrahedron is $\leq V_{oct}$. 

\bibliographystyle{amsplain}

\end{document}